\documentclass[12pt,a4paper,reqno]{article}

\usepackage{amsmath}
\usepackage{amsfonts}
\usepackage{amssymb}
\usepackage{amsthm}
\usepackage{accents}
\usepackage{bbm}
\usepackage{bm}
\usepackage{enumerate}
\usepackage{geometry}
\usepackage{graphicx}
\usepackage[unicode,bookmarks,colorlinks=true,breaklinks=true]{hyperref}
\usepackage[bbgreekl]{mathbbol}
\usepackage[l2tabu, orthodox]{nag}
\usepackage[numbers,square,compress]{natbib}
\usepackage{pgf}
\usepackage{pgfplots}
\usepackage[subrefformat=parens,labelformat=parens]{subfig}
\usepackage{tikz}
\usetikzlibrary{arrows}
\usetikzlibrary{matrix}

\newcommand*{\doi}[1]{\href{http://dx.doi.org/\detokenize{#1}}{doi: \detokenize{#1}}}

\makeatletter
\@ifundefined{chapter}
   {\newtheorem{theorem}{Theorem}[section]}
   {}
\makeatother

\newtheorem*{exercise*}{Exercise}

\DeclareMathSymbol{\widehatsym}{\mathord}{largesymbols}{"62}
\newcommand\lowerwidehatsym{%
  \text{\smash{\raisebox{-1.3ex}{%
    $\widehatsym$}}}}
\newcommand\fixwidehat[1]{%
  \mathchoice
    {\accentset{\displaystyle\lowerwidehatsym}{#1}}
    {\accentset{\textstyle\lowerwidehatsym}{#1}}
    {\accentset{\scriptstyle\lowerwidehatsym}{#1}}
    {\accentset{\scriptscriptstyle\lowerwidehatsym}{#1}}
}

\AtBeginDocument{

    \let\ohat\hat

    \let\bar\undefined
    \let\hat\undefined

    \newcommand{\bar}[1]{\ensuremath{\overline{#1}}}
    \newcommand{\hat}[1]{\ensuremath{\fixwidehat{#1}}}

}

\def\ba{\begin{align}}
\def\ea{\end{align}}

\def\de{\delta}

\newcommand{\mcal}[1]{\ensuremath{\mathcal{#1}}}

\newcommand{\mrm}[1]{\ensuremath{\mathrm{#1}}}

\newcommand{\rsp}{\ensuremath{\mathbb{R}}}

\newcommand{\dd}{\ensuremath{\mathrm{d}}}

\newcommand{\eps}{\ensuremath{\epsilon}}
\newcommand{\veps}{\ensuremath{\varepsilon}}
\newcommand{\phy}{\ensuremath{\varphi}}

\newcommand{\abs}[1]{\ensuremath{\left \vert #1 \right \vert}}

\newcommand{\pair}[2]{\ensuremath{\left < #1 \, , \, #2 \right >}}

\let\div\undefined

\DeclareMathOperator{\curl}{curl}
\DeclareMathOperator{\div}{div}

\newcommand{\unity}{\ensuremath{\mathbb{1}}}

\DeclareMathAlphabet{\mathpzc}{OT1}{pzc}{m}{it}

\newcommand{\functional}[1]{\ensuremath{\mathcal{#1}}}
\newcommand{\matrixform}[1]{\ensuremath{\mathbb{#1}}}

\newcommand{\fa}{\functional{A}}
\newcommand{\fb}{\functional{B}}
\newcommand{\fc}{\functional{C}}
\newcommand{\fg}{\functional{G}}
\newcommand{\fm}{\functional{M}}
\newcommand{\fuu}{\functional{U}}
\newcommand{\cas}{\functional{C}}
\newcommand{\ham}{\functional{H}}
\newcommand{\ffe}{\functional{F}}
\newcommand{\fen}{\functional{S}}
\newcommand{\mom}{\functional{P}}
\newcommand{\nrg}{\functional{E}}

\newcommand{\BJ}{\functional{J}}
\newcommand{\BG}{\functional{G}}

\newcommand{\mass}{\matrixform{M}}
\newcommand{\MA}{\matrixform{A}}
\newcommand{\MG}{\matrixform{G}}
\newcommand{\MI}{\matrixform{I}}

\newcommand{\ML}{\matrixform{L}}
\newcommand{\MS}{\matrixform{S}}
\newcommand{\MC}{\matrixform{C}}

\newcommand{\dist}{\ensuremath{\mathpzc{f}}}
\newcommand{\f}{\dist}
\newcommand{\g}{\ensuremath{\mathpzc{g}}}

\newcommand{\dx}{\ensuremath{\,\dd\xx}}
\newcommand{\dv}{\ensuremath{\,\dd\vv}}
\newcommand{\dz}{\ensuremath{\,\dd\zz}}

\newcommand{\xx}{\ensuremath{x}}
\newcommand{\vv}{\ensuremath{v}}
\newcommand{\zz}{\ensuremath{z}}

\newcommand{\cbra}[2]{\ensuremath{\left \{ #1 , #2 \right \}}}
\newcommand{\sbra}[2]{\ensuremath{\left  [ #1 , #2 \right  ]}}

\newcommand{\fracd}[2]{\ensuremath{\frac{\delta   #1}{\delta   #2}}}
\newcommand{\fracp}[2]{\ensuremath{\frac{\partial #1}{\partial #2}}}

\title{Metriplectic Integrators for\\the Landau Collision Operator}

\author{
\large{Michael Kraus}\\
\small{(michael.kraus@ipp.mpg.de)}
\vspace{.5em}\\
\small{Max-Planck-Institut f\"ur Plasmaphysik}\\
\small{Boltzmannstra\ss{}e 2, 85748 Garching, Deutschland}%
\vspace{.5em}\\
\small{Technische Universit\"at M\"unchen, Zentrum Mathematik}\\
\small{Boltzmannstra\ss{}e 3, 85748 Garching, Deutschland}%
\vspace{1em}\\
\large{Eero Hirvijoki}\\
\small{(ehirvijo@pppl.gov)}
\vspace{.5em}\\
\small{Princeton Plasma Physics Laboratory}\\
\small{Princeton, New Jersey 08543, USA}%
\vspace{1em}\\
}

\date{August 1, 2017}

\begin{document}

\maketitle

\begin{abstract}
  We present a novel framework for addressing the nonlinear Landau
  collision integral in terms of finite element and other subspace
  projection methods. We employ the underlying metriplectic structure
  of the Landau collision integral and, using a Galerkin
  discretization for the velocity space, we transform the
  infinite-dimensional system into a finite-dimensional,
  time-continuous metriplectic system. Temporal discretization is
  accomplished using the concept of discrete gradients. The
  conservation of energy, momentum, and particle densities, as well as
  the production of entropy is demonstrated algebraically for the
  fully discrete system. Due to the generality of our approach, the
  conservation properties and the monotonic behavior of entropy are
  guaranteed for finite element discretizations in general,
  independently of the mesh configuration.
\end{abstract}

\thispagestyle{empty}

\newpage

\tableofcontents

\newpage

\section{Introduction}
Any simulation addressing the Vlasov--Maxwell--Landau system, that
is the system that consists of the Vlasov equation, the Maxwell
equations and the Landau collision integral~\cite{Landau:1936, Landau:1937, Lenard:1960gm}, should respect the
basic laws of physics: total energy should remain constant and entropy a
monotonic function in time, regardless of the total lapse of the
simulation. While the existence of these properties in the continuous
equations is straight forward to demonstrate, no algorithm, at the
time of writing, exists that would retain the same fundamental
properties after spatial and temporal discretizations. Constructing
such an algorithm is one of the outstanding
problems in the topic of numerical simulation of plasmas.

If Coulomb collisions are neglected altogether, recent work on
Vlasov--Poisson, Vlasov--Maxwell and related systems has provided algorithms that satisfy energy
conservation and also preserve other invariants present in the system,
such as the momentum and charge conservation, and the divergence-free
nature of the magnetic field. 
For algorithms in the particle-in-cell framework see for example~\cite{Squire:2012, Evstatiev:2013, Shadwick:2014, Stamm:2014, Xiao:2015, Qin:2016, He:2016, Markidis:2011ep, Chen:2011jf, Chacon:2013cz, Chen:2014eh, Chacon:2016gi, Lapenta:2017uy},
for discontinuous Galerkin methods see~\cite{deDios:2011, deDios:2012a, deDios:2012b, HeathGamba:2012, ChengGamba:2013, Cheng:2014b, Cheng:2014c, Cheng:2014a, Madaule:2014}
and for other grid-based methods see~\cite{Holloway:1996, Back:2014, Kraus:2013:thesis, KrausMajSonnendruecker:2015}.
Realistic kinetic simulations of
plasmas, expanding to macroscopic time scales, however, require also
the effects of Coulomb collisions. The purpose of this paper is to
deliver an algorithm for addressing this issue. We prove algebraically
that our algorithm for discretizing the Landau collision integral in
both space and time conserves the energy, the momentum, and the
particle densities exactly, as well as retains the entropy a
monotonic function in time.

While we base our algorithm on an important recent observation, namely
that the discretization of the weak formulation of the Landau
collision integral with a finite element or any other
subspace-projection method delivers exact conservation laws for
energy, momentum, and particle densities~\cite{Hirvijoki:2017ei}, the
new approach is very different from any previous work on the
subject. We employ the less-familiar representation of the Landau
collision integral in terms of the so-called metriplectic
framework. Similarly as dissipationless Hamiltonian systems can be
described in terms of an energy functional and an antisymmetric
Poisson bracket, the dissipation in metriplectic systems can be
described in terms of an entropy functional and a symmetric, metric
bracket. It is the existence of this framework that ultimately
facilitates the so-called metriplectic discretization of the Landau
collision integral, and allows us to construct an algorithm with the
desired numerical properties.

Metriplectic dynamics is not a commonly encountered topic in the
realms of plasma physics. The first occurrence of a metric bracket for
a single species Landau collision integral, though, dates back to the
pioneering work of Morrison~\cite{Morrison:1986vw}. To bridge the gap,
we start, in Sec.~\ref{sec:metriplectic-dynamics}, by reviewing the
concepts of the metriplectic framework, and outlining a derivation of a
metric bracket for the multispecies Landau operator.
The metric formulation of the collision integral is then
used in conjunction with the more familiar Hamiltonian description of
the Vlasov--Maxwell system to layout a metriplectic formulation for
the kinetic plasma theory. Our discretization of the metric bracket is
discussed in Sec.~\ref{sec:spatial-discretization} and, although
given only for the single-species operator, is extensible to multiple
species in a straight forward way. The discretization employs a
generic finite element method with arbitrary meshing of the velocity
space and is likely extensible to other subspace projection techniques
such as the discontinuous Galerkin approach. After the spatial
discretization is demonstrated, we suggest a temporal discretization
based on the concept of discrete gradients and prove algebraically the
numerical conservation laws and the H-theorem. Numerical
implementation and the extension to discontinuous Galerkin
discretizations is left to be addressed in a future publication.

\section{Metriplectic Dynamics}
\label{sec:metriplectic-dynamics}
Metriplectic dynamics \cite{Morrison:1986vw, Kaufman:1982fl, Kaufman:1984fb, Morrison:1984ca, Morrison:1984wu, Grmela:1984dn, Grmela:1984ea, Grmela:1985jd} provides a convenient framework
for the description of systems that display both Hamiltonian and
dissipative dynamics, such as the Vlasov--Maxwell--Landau system. The
Hamiltonian evolution of the system is determined by a Poisson bracket
$\{ \cdot , \cdot \}$ and the Hamiltonian functional $\ham$, usually
the total energy of the system, and the dissipative evolution is
determined by a metric bracket $( \cdot , \cdot )$ and some functional
$\fen$ that is to be monotonic in time, usually entropy or a generalization
thereof.

In this section, we present the metriplectic formulation of the
Vlasov--Maxwell--Landau system. We first address the general framework, then derive a metric bracket for the multi-species Landau collision integral, and discuss its generalizations. Finally, the metric bracket is used in conjunction
with the Hamiltonian description of the Vlasov--Maxwell system to
provide a metriplectic framework for the kinetic theory in high
temperature plasmas.

\subsection{General framework}
\label{sec:general}
Let us denote by $u (t,z) = (u^{1}, \, u^{2}, \, ..., \, u^{m})^{T}$
the field variables, defined over the domain $\Omega$ with coordinates
$z$, and let $\fuu$ be an arbitrary functional of the field
variables. The domain $\Omega$ is the tensor product of position space
$\Omega_{x}$ and velocity space $\Omega_{v}$. While $\Omega_{x}$ is usually
some bounded domain, possibly periodic in one or more dimensions,
the velocity space $\Omega_{v} = \rsp^{d}$ with $d = \dim \Omega_{x}$.
The evolution of $\fuu$ is given by
\begin{align}
\dfrac{d\fuu}{dt} = \{\fuu,\ffe\} + (\fuu,\ffe) , 
\end{align}
with $\ffe = \ham - \fen$ a generalized free energy functional, analogous with the Gibb's free energy from thermodynamics,  
$\{ \cdot , \cdot \}$ a Poisson bracket, and $( \cdot , \cdot )$ a
metric bracket.  The Poisson bracket, describing the Hamiltonian
evolution, is a bilinear, anti-symmetric bracket of the form
\begin{align}
\cbra{\fa}{\fb} = \int \limits_{\Omega} \dfrac{\de \fa}{\de u^{i}} \, \BJ^{ij} (u) \, \dfrac{\de \fb}{\de u^{j}} \dz ,
\end{align}
where $\fa$ and $\fb$ are functionals of $u$ and $\de \fa / \de u^{i}$
is the functional derivative, defined by
\begin{align}
\dfrac{d}{d\eps} \fa \big[ u^{1}, \, ..., \, u^{i} + \eps v^{i}, \, ..., u^{m} \big] \Big\vert_{\eps=0} = \int \limits_{\Omega} \dfrac{\de \fa}{\de u^{i}} \, v^{i} \dz .
\end{align}
The kernel of the bracket, $\BJ(u)$, is an anti-self-adjoint operator,
which has the property that
\begin{align}
\sum \limits_{l=1}^{m} \left( 
     \frac{\partial \BJ^{ij} (u)}{\partial u^{l}} \, \BJ^{lk} (u)
   + \frac{\partial \BJ^{jk} (u)}{\partial u^{l}} \, \BJ^{li} (u)
   + \frac{\partial \BJ^{ki} (u)}{\partial u^{l}} \, \BJ^{lj} (u)
\right) &= 0 &
& \text{for} &
& 1 \leq i, j, k \leq m ,
\end{align}
ensuring that the bracket $\{ \cdot , \cdot \}$ satisfies the Jacobi
identity,
\begin{align}\label{eq:jacobi_identity}
  \{ \{ \fa , \fb \} , \fc \}
+ \{ \{ \fb , \fc \} , \fa \}
+ \{ \{ \fc , \fa \} , \fb \}
= 0 ,
\end{align}
for arbitrary functionals $\fa , \fb , \fc$ of $u$.  Apart from that,
$\BJ(u)$ is not required to be of any particular form, and most
importantly it is allowed to depend on the fields $u$.  If $\BJ(u)$ has
a non-empty nullspace, there exist so-called Casimir invariants, that
is functionals $\cas$ for which $\{ \fa , \cas \} = 0$ for all
functionals $\fa$. The monotonic entropy functional $\fen$ is usually
one of these Casimir invariants.

The metric bracket $( \cdot , \cdot )$, describing dissipative
effects, is a symmetric bracket, defined in a similar way as the
Poisson bracket by
\begin{align}
( \fa , \fb ) = \int \limits_{\Omega} \dfrac{\de \fa}{\de u^{i}} \, \BG^{ij} (u) \, \dfrac{\de \fb}{\de u^{j}} \dz ,
\end{align}
where $\BG(u)$ is now a self-adjoint operator with an appropriate nullspace
such that $( \ham , \ffe ) = 0$. All Casimirs $\cas$ of the Poisson
bracket should be Casimirs also of the metric bracket, except the
functional $\fen$ which is explicitly required not to be a Casimir of the metric bracket.

In this paper, we choose a convention that dissipates the free energy, conserves the Hamiltonian, and produces entropy. Other conventions are equally possible, but we stick with the concepts familiar from thermodynamics, so that the equilibrium state is reached when free energy is at minimum. For this framework and our conventions to be consistent, it is essential that (i) $\ham$ is a Casimir of the metric bracket, (ii) $\fen$ is a Casimir of the Poisson bracket, and that (iii) the metric bracket is negative semi-definite. With respect to these choices, we then have
\begin{align}
\dfrac{d\ham}{dt} &= \{\ham,\ffe\}+(\ham,\ffe)=\{\ham,-\fen\}=0,\\
\dfrac{d\fen}{dt} &= \{\fen,\ffe\}+(\fen,\ffe)=-(\fen,\fen)\ge 0,\\
\dfrac{d\ffe}{dt} &= \{\ffe,\ffe\}+(\ffe,\ffe)=(\ffe,\ffe)\le 0,
\end{align}
reproducing the First and Second Laws of Thermodynamics and dissipation of free energy.

For an equilibrium state $u_{eq}$, the time evolution of any functional $\fuu$ stalls, the free-energy functional reaches its minimum, and entropy production stops. If the metriplectic system has no Casimirs, the equilibrium state satisfies an energy principle, according to which the first variation vanishes, $\delta \ffe[u_{eq}]=0$, and the second variation is strictly positive, $\delta^2\ffe[u_{eq}]> 0$ (for details see e.g.~\citet{Holm:1985}).
If Casimirs $\fc_i$ exist, the equilibrium state becomes degenerate, and the energy principle must be modified to account for the existing Casimirs. This leads to the so-called energy--Casimir principle~\cite{Morrison:1998}. In this case the the equilibrium state satisfies
\begin{align}\label{eq:energy-casimir-principle}
\delta \ffe[u_{eq}]+\sum_i\lambda_i\delta \fc_i[u_{eq}]=0,
\end{align}
where $\lambda_i$ act as Lagrange multipliers and are determined uniquely from the mass, momentum and energy of the initial conditions for $u$.
Lastly, for each $z \in \Omega$ the
equilibrium state of the free-energy functional $\ffe$ is unique.
This can be accomplished by employing a convexity argument, namely
if $\Omega$ is a convex domain and $\ffe$ is strictly convex,
then $\ffe$ has at most one critical point (see e.g. \citet{GiaquintaHildebrandt:2004} for details).
This is the case if 
\begin{align}
\delta^{2} ( \ffe + \sum_{i} \lambda_{i} \cas_{i} ) > 0,
\end{align}
for the non-vanishing field $u_{eq}$.

\subsection{Metric bracket for Multi-species Landau Collision Integral}
\label{sec:landau_multi_species}

In high-temperature plasmas, collisional evolution of the distribution
function $\f_s$ of species $s$ due to collisions with multiple other
species $s'$ (including $s$) can be described by the Landau collision
integral
\begin{equation}\label{eq:landau-operator_multi-species}
\frac{\partial \f_s}{\partial t} = C[\f_s](v)=\sum_{s'}\frac{c_{ss'}}{m_s}\frac{\partial}{\partial v}\cdot \int \limits_{\Omega_{v}} U(v,v') \cdot J_{ss'}(v,v') \dv',
\end{equation}
which is a nonlinear integro-differential operator in velocity
space. The antisymmetric vector $J_{ss'}(v,v')=-J_{s's}(v',v)$ in the
operator depends on $f_{s}$ and $f_{s'}$ and is defined as
\begin{equation}\label{eq:landau_Jvector}
J_{ss'}(v,v')
= \frac{\f_{s'}(v')}{m_{s }}  \fracp{\f_{s }(v )}{v }
- \frac{\f_{s }(v )}{m_{s'}}  \fracp{\f_{s'}(v')}{v'} ,
\end{equation}
and the symmetric parameter $c_{ss'}$ is given by
\begin{equation}
c_{ss'}=\frac{q_s^2q_{s'}^2\ln \Lambda}{8\pi\varepsilon_0^2},
\end{equation}
with $q_s$ and $m_s$ the charge and mass of particles of species $s$, $\varepsilon_0$ the vacuum permittivity, 
and $\ln\Lambda$ the Coulomb logarithm.
The Landau tensor $U (v, v')$,
valid at non-relativistic energies, is a scaled projection matrix of
the relative velocity $v - v'$ between the colliding particles,
\begin{align}\label{eq:landau-tensor}
 U_{ij} (v, v') = \dfrac{1}{\abs{v - v'}} \left( \delta_{ij} - \dfrac{(v_{i} - v_{i}') (v_{j} - v_{j}')}{\abs{v - v'}^2} \right) .
\end{align}
Let us multiply with a function $\g_s$ and integrate over the infinite
velocity space. Carrying out a partial integration provides us with
the expression
\begin{equation}\label{eq:landau_weak}
\int \limits_{\Omega_{v}}\g_s(v)\frac{\partial \f_s}{\partial t} \dv=-\sum_{s'}\frac{c_{ss'}}{m_s}\int \limits_{\Omega_{v}}\frac{\partial \g_s(v)}{\partial v}\cdot \int \limits_{\Omega_{v}} U(v,v')\cdot J_{ss'}(v,v') \dv'\dv .
\end{equation}
Exchanging the species indices, we may similarly write for another species
\begin{align}
\int \limits_{\Omega_{v}}\g_{s'}(v)\frac{\partial \f_{s'}}{\partial t} \dv
&=\sum_{s}\frac{c_{ss'}}{m_{s'}}\int \limits_{\Omega_{v}}\int \limits_{\Omega_{v}}\frac{\partial \g_{s'}(v')}{\partial v'}\cdot  U(v,v')\cdot J_{ss'}(v,v')\dv'\dv ,
\end{align}
where the symmetry of $c_{ss'}$ and $U(v,v')$ and the
antisymmetry of $J_{ss'}(v,v')$ were used. If we then sum over the
different species, we obtain
\begin{align}\label{eq:expectation-sum}
\sum_{s}\int \limits_{\Omega_{v}}\g_{s}(v)\frac{\partial \f_{s}}{\partial t} \dv & = -\sum_{s,s'}\frac{c_{ss'}}{2}\int \limits_{\Omega_{v}}\int \limits_{\Omega_{v}}\left(\frac{1}{m_{s}}\frac{\partial \g_{s}(v)}{\partial v}-\frac{1}{m_{s'}}\frac{\partial \g_{s'}(v')}{\partial v'}\right)\nonumber\\& \qquad \cdot  U(v,v')\cdot J_{ss'}(v,v')\dv'\dv.
\end{align}
It is straightforward to verify that if $\g_{s}(v)=m_s\{1,v,|v|^2\}$
the integral on the right-hand side vanishes identically. Since the left-hand side
then presents the total collisional rate of change in time of mass
density, kinetic momentum density, and kinetic energy density in the
plasma, the above calculation verifies that these three quantities are
invariants of the Landau collision operator.

Starting from the weak formulation, that is Equation~\eqref{eq:landau_weak} with $\g_s$ interpreted as a test function, the steps to derive a metric
bracket for the collision integral are straightforward. We introduce
the functionals
\begin{align}
\label{eq:entropy-vspace}
\fen&=-T\sum_{s}\int \limits_{\Omega_{v}} \f_s(t,v) \ln \bigg( \dfrac{\f_s(t,v)}{g} \bigg) \dv , \\
\label{eq:g-functional-vspace}
\fg&=\sum_s \fg_s, \qquad
\fg_s=\int \limits_{\Omega_{v}} \g_s(v) \, \f_s(t,v) \dv,
\end{align}
where $T$ is a normalization constant in units of energy density and $g$ is a normalization factor for the distribution function. The above functionals allow us to formally rewrite the vector $J_{ss'}(v,v')$ according to
\begin{align}
J_{ss'}(v,v')=-\frac{\f_{s'}(v')\f_{s}(v)}{T}\left(\frac{1}{m_{s }} \fracp{}{v}\frac{\delta \fen }{\delta \f_s} 
- \frac{1}{m_{s'}}\fracp{}{v'}\frac{\delta \fen}{\delta \f_{s'}}\right).
\end{align}
Following the requirement for the metric bracket to be negative semi-definite, the above definitions allow us to re-express \eqref{eq:expectation-sum} in the form
\begin{align}
\label{eq:landau-evolution}
\frac{\partial \fg}{\partial t}&=(\fg,-\fen),
\end{align}
where the symmetric bracket is defined according to
\begin{align}\label{eq:simple-symmetric-bracket}
(\fa,\fb)& =-\sum_{s,s'}\frac{c_{ss'}}{2T}\int \limits_{\Omega_{v}}\int \limits_{\Omega_{v}}\left(\frac{1}{m_{s }}  \fracp{}{v }\frac{\delta \fa }{\delta \f_s} 
- \frac{1}{m_{s'}}\fracp{}{v'}\frac{\delta \fa }{\delta \f_{s'}}\right)\nonumber\\&\qquad\cdot  U(v,v')\f_{s'}(v')\f_{s }(v )\cdot\left(\frac{1}{m_{s }}  \fracp{}{v }\frac{\delta \fb }{\delta \f_s} 
- \frac{1}{m_{s'}}\fracp{}{v'}\frac{\delta \fb }{\delta \f_{s'}}\right)\dv'\dv,
\end{align}
and is readily identified as the metric bracket for the Landau collision integral (cf. Ref.~\cite{Morrison:1986vw}). 
It is straightforward to verify that inserting~\eqref{eq:entropy-vspace} and~\eqref{eq:g-functional-vspace} into~\eqref{eq:landau-evolution} leads to a weak formulation of equation~\eqref{eq:landau-operator_multi-species}. Further, defining the mass, momentum and kinetic energy functionals by
\begin{align}
\fm_s &= m_s\int \limits_{\Omega_{v}}\,\f_s(t,v) \dv, &
\mom &= \sum_s m_s\int \limits_{\Omega_{v}}\, v\, \f_s(t,v) \dv, &
\nrg &= \sum_s \frac{m_s}{2}\int \limits_{\Omega_{v}} \, \vert v \vert^2 \, \f_s(t,v) \dv,
\end{align}
it is straightforward to verify that the bracket satisfies
\begin{align}
(\fm_s,\fa)=0, \qquad
({\cal P},\fa)=0, \qquad
(\nrg,\fa)=0,
\end{align}
for arbitrary functionals $\fa$. Thus $\fm_s$, $\mom$, and $\nrg$ are Casimirs of the bracket~\eqref{eq:simple-symmetric-bracket}. Restricting our attention only to the velocity space, for now, we may then conclude that the collisional dynamics of a given functional $\fa$, consistent with the Landau collision integral, follow from 
\begin{align}
\dfrac{d\fa}{dt}=(\fa,\ffe),
\end{align}
where the free-energy functional is $\ffe=\nrg-\fen$. This formulation satisfies the basic principles of thermodynamics, i.e., $d\ffe /dt\le 0$, $d\nrg/dt=0$, and $d\fen/dt\ge 0$. 
The equilibrium state obeys the energy-Casimir principle~\eqref{eq:energy-casimir-principle}, which in this case provides
\begin{align}
\bigg(\frac{\delta \fen}{\delta \f_s}+\lambda_{s}\frac{\delta \fm_s}{\delta \f_s}+\lambda_{{\cal P}}\cdot\frac{\delta \mom}{\delta \f_s}+\lambda_{{\cal E}}\frac{\delta \nrg}{\delta \f_s} \bigg) \bigg\vert_{\f = \f_{eq}}=0, \quad \text{for all} \quad s,
\end{align}
and leads to the following condition for the equilibrium distribution functions
\begin{align}
-T \bigg( 1 + \ln \bigg( \dfrac{\f_s(t,v)}{g} \bigg) \bigg) + \lambda_s m_s + \lambda_{\cal P}\cdot m_s v + \lambda_{\cal E}\frac{m_s}{2}\vert v\vert^2 = 0, \quad \text{for all} \quad s.
\end{align}
The equilibrium distributions with respect to collisional dynamics are thus identified as Maxwellians with each species having common temperature and flow velocity but possibly different densities.

\subsection{Metric Bracket for General Collision Operators}

The fact that the single-species Landau collision operator, as well as
other small-angle Coulomb collision operators relevant for plasmas,
can be obtained from a general metric bracket was demonstrated
by Morrison already in 1986~\cite{Morrison:1986vw}. Morrison's bracket can
easily be generalized to the multi-species case,
\begin{align}\label{eq:metric_bracket_collision_operator}
( \fa , \fb )[\f] = \sum_{s',s''} \int \limits_{\Omega} \int \limits_{\Omega} \Gamma_{s's''}(\fa;z',z'') \cdot T_{s's''}(z'; z'') 
\cdot \Gamma_{s's''}(\fb;z',z'') \dz' \dz'' ,
\end{align}
where the three-component vector $\Gamma_{s's''}(\fa;z',z'')$ is defined as
\begin{align}
\Gamma_{s's''}(\fa;z',z'')=\dfrac{1}{m_{s'}}\dfrac{\partial}{\partial v'} \dfrac{\de \fa}{\de \f_{s'} (z')} - \dfrac{1}{m_{s''}}\dfrac{\partial}{\partial v''} \dfrac{\de \fa}{\de \f_{s''} (z'')} ,
\end{align}
and $T_{s's''}(z';z'')=W_{s's''}(z';z'')\delta(x'-x'')$ with $W_{s's''}$ a symmetric,
positive or negative semi-definite matrix with an eigenvector $v'-v''$
corresponding to a zero eigenvalue. The choice between positive and negative depends on the choice of dissipating either the free energy or entropy.
Different collision operators follow from different choices for the matrix $W_{s's''}$ and the entropy
functional $\fen$, which is restricted to be of the form
\begin{align}\label{eq:metric_bracket_entropy_functional}
\fen = \sum \limits_{s} \int \limits_{\Omega} s(\f_{s}) \dz ,
\end{align}
with $s$ an arbitrary function of $\f$, and required to be a Casimir of the Poisson bracket.

For kinetic systems, such as the Vlasov--Maxwell system, the total energy and momentum are typically given in terms of the kinetic energy and momentum of the plasma particles and some field energy $\ham_{\mrm{fields}}$ and momentum terms $\mom_{\mrm{fields}}$ according to 
\begin{align}\label{eq:kinetic_hamiltonian}
\ham[\f] &= \sum\limits_{s}\dfrac{m_s}{2} \int \limits_{\Omega} \, \vert v\vert^2 \, \f_s(t,z) \dz + \ham_{\mrm{fields}},\\
\label{eq:kinetic_momentum}
\mom[\f] & = \sum_s m_s\int \limits_{\Omega} \, v \, \f_s(t,z)\, dz + \mom_{\mrm{fields}}.
\end{align}
Even if the field-energy and field-momentum were not independent of $\f$, the total Hamiltonian $\ham$ and momentum $\mom$ will be Casimirs of the symmetric bracket as long as their functional derivatives satisfy the conditions
\begin{align}
\frac{\delta \ham}{\delta \f_s}=\frac{m_s}{2}\vert v\vert^2, \qquad \frac{\delta \mom}{\delta \f_s}=m_s\,v.
\end{align}
Similarly, a generalized mass will be a Casimir of the metric bracket as long as $\delta \fm / \delta \f_s=m_s$. These results follow directly from the conditions 
\begin{align}
\Gamma_{s's''}(\fm;z',z'')=0,\quad \Gamma_{s's''}(\mom;z',z'')=0, \quad \Gamma_{s's''}(\ham;z',z'')\cdot W_{s's''}(z',z'')=0.
\end{align}
One may thus define a free-energy functional $\ffe=\ham+\lambda_{\fen}\, \fen$, and, to conclude, that the collisional evolution of any functional $\fuu$ is determined by
\begin{align}
\left(\frac{\partial \fuu}{\partial t}\right)_c=(\fuu,\ffe).
\end{align}
Maxwell--Boltzmann statistics and the Landau collision integral, for example, are obtained by choosing $\lambda_{\fen}=-1$, $s(\f)=-T\f\ln\f$, and
\begin{align}\label{eq:metriplectic_W}
W_{s's'',ij}(z',z'')= -\frac{c_{s' s''}}{2T} U_{ij} (v', v'') \, M(\f_{s'}(z'))\,M(\f_{s''}(z'')),
\end{align}
with $M(\f)=\f$, while Fermi-Dirac statistics and the relevant collision operator follow from choosing $s(\f)=-T(\f \ln \f+(1-\f)\ln(1-\f))$ and $M(\f)=\f(1-\f)$, as discussed by~\citet{Morrison:1986vw}.
Following the energy-Casimir principle, the general equilibrium condition is obtained from
\begin{align}
\bigg(\frac{\delta \ham}{\delta \f_s} +\lambda_{\fen} \frac{\delta \fen}{\delta \f_s} + \lambda_{\ham}\frac{\delta\ham}{\delta \f_s}+\lambda_{\mom}\frac{\delta \mom}{\delta \f_s}+\lambda_{s}\frac{\delta \fm_s}{\delta \f_s}\bigg) \bigg\vert_{\f = \f_{eq}}=0
\end{align}
which leads to
\begin{align}
\lambda_{\fen}s_{\f}(\f_{s,eq})+(\lambda_{\ham}+1)\frac{m_s}{2}\vert v\vert^2+\lambda_{\mom}\cdot m_sv+\lambda_{s}m_s=0
\end{align}
Since the second variations of $\ham$, $\mom$, and $\fm_s$ with respect to $\f_s$ all vanish, the equilibrium is unique as long as $\delta^2\fen < 0$. This is certainly true for both the Maxwell--Boltzmann and Fermi--Dirac entropies.

It is worth pointing out that these results are fully general and hold
for all choices of $W_{s's''}$, which satisfy the assumptions that $W_{s's''}$ is
symmetric, positive or negative semi-definite, depending on the convention, and have a zero eigenvector $v'-v''$, i.e.,
\begin{subequations}\label{eq:metric-bracket-symmetries}
\begin{align}
W_{s's'',ij} (z', z'') &= W_{s's'',ji} (z', z'') , \\
W_{s's'',ij} (z', z'') &= W_{s's'',ij} (z'', z') , \\
(v_{i}' - v_{i}'') \, W_{s's'',ij} &= 0 .
\end{align}
\end{subequations}
Hence the equilibrium state only depends
on the choice of the entropy functional $\fen$ but not on the choice
of $W_{s's''}$.  In particular, it turns out that there is no strict need for
a compatibility condition
between $M$ in~\eqref{eq:metriplectic_W} and the function
$s$ in~\eqref{eq:metric_bracket_entropy_functional} or its derivatives
as mentioned in~\cite{Morrison:1986vw, Morrison:2017}.  
However, the freedom to choose $M$ often has to be used to regularize the
collision operator. One often finds the distribution function in the
denominator of the second derivative of the entropy function $s_{\f\f}$.
In order to avoid singularities, $M$ can be chosen appropriately to
cancel such denominators.
In other cases, it might also be possible to choose $M$ in order to
simplify the bracket, e.g., to reduce the degree of the nonlinearity
in the distribution function $\f$ and therefore improve the
computational tractability.

\subsection{Metriplectic Vlasov--Maxwell--Landau System}

The Vlasov--Maxwell--Landau system of equations is given by the set
\begingroup
\allowdisplaybreaks
\begin{subequations}\label{eq:vlasov-maxwell}
\begin{align}
\fracp{\f_s}{t} + \vv \cdot \nabla_\xx \f_s + \frac{q_s}{m_s} (E + \vv \times B) \cdot \nabla_\vv \f_s &= C[\f_s](v), \label{eq:vlasov}
\\
\noalign{\medskip}
\frac{1}{c^2}\fracp{E}{t} - \curl B &= - \mu_0j, \label{eq:ampere}
\\
\noalign{\medskip}
\fracp{B}{t} + \curl E &= 0, \label{eq:faraday}
\\
\noalign{\medskip}
\div E &= \frac{\rho}{\varepsilon_0}, \label{eq:gauss}
\\
\noalign{\medskip}
\div B &= 0, \label{eq:gaussm}
\end{align}
\end{subequations}
\endgroup
where $c=\sqrt{\varepsilon_0\mu_0}^{-1}$ is the speed of light, $\mu_0$ the vacuum permeability, $\varepsilon_{0}$ is the vacuum permittivity, $E$ and $B$ denote the electric and magnetic fields, $\rho$ is
the charge density and $j$ the current density, defined in terms of
the distribution functions $\f_s$ according to
\begin{align}\label{eq:sources}
\rho &= \sum_s q_s \int \limits_{\Omega_{v}} \f_s \dv, &
j &= \sum_s q_s \int \limits_{\Omega_{v}} \vv \f_s \dv.
\end{align}
This system of equations can be obtained alternatively from a
metriplectic system as
\begin{align}
\frac{\partial \fuu}{\partial t} = \cbra{\fuu}{\ffe}+(\fuu,\ffe),
\end{align}
described by a noncanonical Poisson bracket~\cite{Morrison:1980,
  Weinstein:1981, Marsden:1982, Morrison:1982, Morrison:2013},
\begin{align}\label{eq:poisson_bracket_maxwell}
\cbra{\fa}{\fb} [\f, E, B]
\nonumber
=& \sum \limits_{s} \int \limits_{\Omega_x} \int \limits_{\Omega_v} \dfrac{\f_s}{m_s} \sbra{ \fracd{\fa}{\f_s} }{ \fracd{\fb}{\f_s} } \dx \dv \\
\nonumber
&+ \dfrac{1}{\veps_{0}} \sum \limits_{s} \dfrac{q_s}{m_s} \int \limits_{\Omega_x} \int \limits_{\Omega_v} \f_s \left(
        \nabla_\vv \fracd{\fa}{\f_s} \cdot \fracd{\fb}{E}
      - \nabla_\vv \fracd{\fb}{\f_s} \cdot \fracd{\fa}{E}
   \right) \dx \dv \\
\nonumber
&+ \sum \limits_{s} \dfrac{q_s}{m_s^2} \int \limits_{\Omega_x} \int \limits_{\Omega_v} \f_s \, B \cdot \left( \nabla_\vv \fracd{\fa}{\f_s} \times \nabla_\vv \fracd{G}{\f_s} \right) \dx \dv \\
&+ \dfrac{1}{\veps_{0}} \int \limits_{\Omega_x} \left(
        \curl \fracd{\fa}{E} \cdot \fracd{\fb}{B}
      - \curl \fracd{\fb}{E} \cdot \fracd{\fa}{B}
   \right) \dx ,
\end{align}
a symmetric, negative semi-definite metric bracket
\begin{align}
( \fa , \fb )[\f] = \sum_{s',s''} \int \limits_{\Omega} \int \limits_{\Omega} \Gamma_{s's''}(\fa;z',z'') \cdot T_{s's''}(z'; z'') 
\cdot \Gamma_{s's''}(\fb;z',z'') \dz' \dz'' ,
\end{align}
a Hamiltonian functional $\ham$, given by the sum of the kinetic energy and
the electric and magnetic field energies, namely
\begin{align}
\ham[\f,E,B] = \sum \limits_{s} \dfrac{m_s}{2} \int \limits_{\Omega_x} \int \limits_{\Omega_v} \abs{\vv}^{2} \, \f_s (\xx,\vv) \dx \dv + \dfrac{1}{2} \int \limits_{\Omega_{x}} \Big( \varepsilon_0\abs{E (\xx)}^{2} + \mu_0^{-1}\abs{B (\xx)}^{2} \Big) \dx ,
\label{eq:hamiltonian_vlasov_maxwell}
\end{align}
an entropy functional
\begin{align}
\fen[\f] = -T\sum \limits_s \int \limits_{\Omega} \f_s(z)\,\ln \left(\frac{\f_s(z)}{g}\right) \dz,
\end{align}
where $T$ is the temperature and $g$ a normalization constant, and a free-energy functional 
\begin{align}
\ffe[\f,E,B]=\ham[\f,E,B]-\fen[\f].
\end{align}
The canonical single-particle Poisson bracket, present in the functional Poisson bracket, is
\begin{align}
[f,g]=\nabla_\xx f \cdot  \nabla_\vv g - \nabla_\xx g \cdot  \nabla_\vv f,
\end{align}
and the vector $\Gamma_{s's''}(\fa;z',z'')$ and the matrix $T_{s's''}(z'; z'')$, present in the metric bracket, are
\begin{align}
\Gamma_{s's''}(\fa;z',z'')&=\dfrac{\nabla_{v'}}{m_{s'}}\dfrac{\de \fa}{\de \f_{s'} (z')} - \dfrac{\nabla_{v''}}{m_{s''}} \dfrac{\de \fa}{\de \f_{s''} (z'')},\\
T_{s's''}(z'; z'')&=-\frac{c_{s's''}}{2T}\delta(x'-x'')U(v',v'')\f_{s'}(z')\f_{s''}(z'').
\end{align}
The sign conventions are chosen so that free energy is dissipated, $d\ffe/dt\le 0$, total energy is conserved, $d\ham/dt=0$, and entropy is produced $d\fen/dt\ge 0$.

The total momentum, consisting of the kinetic and field contributions,  
\begin{align}
\label{eq:momentum}
\mom[\f,E,B] & = \sum_s m_s \int \limits_{\Omega_{x}} \int \limits_{\Omega_{v}} v \, \f_s \dx \dv +\mu_0^{-1} \int \limits_{\Omega_{x}} E \times B \dx,
\end{align}
is a Casimir of the metric bracket but not of the Poisson bracket. It is conserved, though, given that the Poisson equation for the electric field is satisfied. The only Casimirs of this full system are the species-wise mass functionals,
\begin{align}
\fm_s=m_s \int \limits_{\Omega} \f_s(z)\, dz,
\end{align}
and the equilibrium state is obtained from the modified energy--Casimir principle,
\begin{align}
\de \ffe + \sum \limits_s \lambda_s \, \de \fm_s
&= \sum_s\left( \dfrac{\de \ham}{\de \f_s} - \dfrac{\de \fen}{\de \f_s} + \lambda_s \dfrac{\de \fm_s}{\de \f_s} \right) \de \f_{s,eq} + \dfrac{\de \ham}{\de E} \cdot \de E_{eq} + \dfrac{\de \ham}{\de B} \cdot \de B_{eq}
= 0. 
\end{align}
For arbitrary variations $\de \f_{s,eq}$, $\de E_{eq}$ and $\de B_{eq}$, the principle leads to
\begin{align}\label{eq:metriplectic_bracket_equilibrium_2}
-T \bigg( 1 + \ln \bigg( \dfrac{\f_{s,eq}}{g} \bigg) \bigg) &= \frac{m_s}{2}\vert v\vert^{2} + \lambda_sm_s , &
E_{eq} (x) &= 0 , &
B_{eq} (x) &= 0 ,
\end{align}
and describes Maxwellian distribution functions with zero flow and equal uniform temperature. The numbers $\lambda_s$ are uniquely determined from the initial conditions by
\begin{align*}
\int \limits_{\Omega} \f_{0,s} (z) \, dz = \int \limits_{\Omega} \f_{s,eq} (z) \, dz .
\end{align*}

\section{Spatial Discretization}
\label{sec:spatial-discretization}
In the following, we will restrict our attention to the spatial and
temporal discretization of the metric part of the system. 
The distribution function $\f (t,v)$ is considered as a function of
time and velocity only. The full Vlasov--Maxwell--Landau system will
be considered in subsequent publications~\cite{Kraus:2017:VI-VML, Kraus:2017:GEMPIC-VMFP}.
We shall also clarify the discussion assuming only one species and using normalized units.
The metric bracket addressed in this section can then be expressed as
\begin{multline}\label{eq:single-species-bracket}
( \fa , \fb )[\f] = -\frac{1}{2}\int \limits_{\Omega} \int \limits_{\Omega} \bigg( \dfrac{\partial}{\partial v'} \dfrac{\de \fa}{\de \f (v')} - \dfrac{\partial}{\partial v''} \dfrac{\de \fa}{\de \f (v'')} \bigg) \cdot M(\f(v')) \, U(v'; v'') \, M(\f(v'')) \\
\cdot \bigg( \dfrac{\partial}{\partial v'} \dfrac{\de \fb}{\de \f (v')} - \dfrac{\partial}{\partial v''} \dfrac{\de \fb}{\de \f (v'')} \bigg) \dv' \dv'' .
\end{multline}
Let us note, though, that the generalisation to multiple species is straight forward.

Our approach for the spatial discretization is a combination of the ones
presented in~\cite{Kraus:2016} and~\cite{Hirvijoki:2017ei}.
We will consider a finite dimensional space
$Q_h(\Omega) \subset L^{2}(\Omega)$, which is a subset of all square
integrable functions, defined on the domain $\Omega$. 
The discrete domain has a tensor product structure, similar to the
continuous domain, but in contrast to the continuous velocity space
the discrete velocity space is often bounded. We denote the basis
functions in this space by $\phy_{i} (v)$, so that an element
$\f_{h} \in Q_h(\Omega)$, approximating the distribution function
$\f \in L^{2}(\Omega)$, can be expressed as
\begin{align}
\f_{h} (t,v) = \sum \limits_{i=1}^{N} \ohat{f}_{i} (t) \, \phy_{i} (v) ,
\end{align}
where the $\ohat{f}_{i} (t)$ denote the degrees of freedom. The
configuration of the system at a given point in time is then determined by
$\ohat{f} (t) = (\ohat{f}_{1} (t), ..., \ohat{f}_{N} (t))^{T} \in
\rsp^{N}$.

\subsection{Discrete Functional Derivatives}

Consider some functional $\fa$ of the distribution function $\f$.  Its
functional derivative with respect to $\f$ is defined by
\begin{align}\label{eq:dfd_continuous_functional_derivative}
\dfrac{d}{d\eps} \fa \big[ \f + \eps \g \big] \Big\vert_{\eps=0}
= \pair{\dfrac{\de \fa}{\de \f}}{\g}_{L^{2}}
= \int \limits_{\Omega} \dfrac{\de \fa}{\de \f} \, \g (v) \dv .
\end{align}
Here, $\g$ is an element of the same space as $\f$, that is
$\g \in L^{2} (\Omega)$, while the functional derivative
$\de \fa / \de \f$ is an element of the dual space of
$L^{2} (\Omega)$, and $\pair{\cdot}{\cdot}$ denotes
the $L^{2}$ inner product.
When we apply the functional $\fa$ to the Galerkin approximation
of $\f$, it becomes a function $\hat{A}$ of the degrees of freedom
$\ohat{f}$,
\begin{align}
\fa [\f_{h}] = \hat{A} (\ohat{f}) .
\end{align}
In order to obtain a discrete version of the metric bracket, the
functional derivative $\de \fa / \de \f$ needs to be replaced with a
partial derivative $\partial \hat{A} / \partial \ohat{f}$. To do so,
we require that the pairing
in~\eqref{eq:dfd_continuous_functional_derivative} be equal to some
finite-dimensional equivalent, that is
\begin{align}\label{eq:dfd_discrete_functional_derivative}
\pair{\dfrac{\de \fa [\f_{h}]}{\de \f}}{\g_{h}}_{L^{2}}
= \pair{\dfrac{\partial \hat{A}}{\partial \ohat{f}}}{\ohat{g}}_{\rsp^{N}} 
= \sum \limits_{i=1}^{N} \dfrac{\partial \hat{A}}{\partial \ohat{f}_{i}} \, \ohat{g}_{i} ,
\end{align}
where $\ohat{g} (t) = (\ohat{g}_{1} (t), ..., \ohat{g}_{N} (t))^{T}$
denotes the degrees of freedom of $\g_{h}$, so that
\begin{align}
\g_{h} (t, v) = \sum \limits_{i=1}^{N} \ohat{g}_{i} (t) \, \phy_{i} (v) .
\end{align}
Let us denote the dual basis to
$\phy (v) = (\phy_{1} (v), ..., \phy_{N} (v))^{T}$ in $L^{2}$ by
$\psi (v) = (\psi_{1} (v), ..., \psi_{N} (v))^{T}$, so that
\begin{align}\label{eq:dfd_dual_basis_delta}
\int \limits_{\Omega} \psi_{i} (v) \, \phy_{j} (v) \dv = \delta_{ij}
\quad \text{for} \quad
1 \leq i, j \leq N .
\end{align}
In the dual basis, the functional derivative
in~\eqref{eq:dfd_discrete_functional_derivative} can be written as
\begin{align}\label{eq:dfd_functional_derivative_dual_basis}
\dfrac{\de \fa [\f_{h}]}{\de \f} = \sum \limits_{i=1}^{N} a_i \, \psi_{i} (v) ,
\end{align}
where the degrees of freedom $a_i$ are still to be determined.
Using \eqref{eq:dfd_discrete_functional_derivative} and
\eqref{eq:dfd_dual_basis_delta} for
$\ohat{g} = (0, \, \dots, \, 0, \, 1, \, 0, \, \dots, \, 0)^\top$ with
$1$ at the $i$-th position and $0$ everywhere else, so that
$\g_{h} = \phy_{i}$, we find that
\begin{align}
a_{i} = \dfrac{\partial \hat{A}}{\partial \ohat{f}_{i}} ,
\end{align}
and can thus write
\begin{align}\label{eq:dfd_functional_derivative_dual_basis2}
\dfrac{\de \fa [\f_{h}]}{\de \f} = \sum \limits_{i=1}^{N} \dfrac{\partial \hat{A}}{\partial \ohat{f}_{i}} \, \psi_{i} (v) .
\end{align}
Expressing the dual basis $\psi$ in terms of the primal basis $\phy$
as
\begin{align}\label{eq:dfd_dual_basis_in_primal_basis}
\psi_{i} (v) = \sum \limits_{j=1}^{N} \MA_{ij} \, \phy_{j} (v) ,
\end{align}
this becomes
\begin{align}\label{eq:dfd_functional_derivative_primal_basis}
\dfrac{\de \fa [\f_{h}]}{\de \f} = \sum \limits_{i,j=1}^{N} \dfrac{\partial \hat{A}}{\partial \ohat{f}_{i}} \, \MA_{ij} \, \phy_{j} (v) .
\end{align}
In order to determine the unknown coefficients $\MA_{ij}$, let us compute the
$L_{2}$ inner product of \eqref{eq:dfd_dual_basis_in_primal_basis}
with $\phy_k$,
\begin{align}
\int \limits_{\Omega} \psi_{i} (v) \, \phy_{k} (v) \dv
= \int \limits_{\Omega} \sum \limits_{j=1}^{N} \MA_{ij} \, \phy_{j} (v) \, \phy_{k} (v) \dv
= \sum \limits_{j=1}^{N} \MA_{ij} \int \limits_{\Omega} \phy_{j} (v) \, \phy_{k} (v) \dv .
\end{align}
Denoting by $\mass$ the mass matrix of the basis functions $\phy$,
\begin{align}
\mass_{jk} = \int \limits_{\Omega} \phy_{j} (v) \, \phy_{k} (v) \dv ,
\end{align}
and using~\eqref{eq:dfd_dual_basis_delta}, we obtain the relation $\MI = \MA \mass$ with $\MI$ the $N \times N$ identity matrix,
so that $\MA = \mass^{-1}$ and therefore
\begin{align}\label{eq:discrete_functional_derivative}
\dfrac{\de \fa [\f_{h}]}{\de \f} = \sum \limits_{i,j=1}^{N} \dfrac{\partial \hat{A}}{\partial \ohat{f}_{i}} \, (\mass^{-1})_{ij} \, \phy_{j} (v) .
\end{align}
This is the relation we have sought for, expressing the ``continuous''
functional derivative $\de \fa / \de \f$ in terms of the ``discrete''
partial derivative $\partial \hat{A} / \partial \ohat{f}$.

\subsection{Discrete Bracket and Semi-discrete Equations of Motion}

In order to obtain a discrete version of the metriplectic
bracket, we restrict~\eqref{eq:single-species-bracket} to the space
of functionals on $Q_h(\Omega)$. This allows us to replace the
functional derivatives with partial derivatives as in~\eqref{eq:discrete_functional_derivative}
and compute the remaining integrals, which leads to
\begin{align}\label{eq:discrete_metric_bracket_collision_operator}
( \hat{A} , \hat{B} )_{h} = \sum \limits_{i,j,k,\ell=1}^{N}\dfrac{\partial \hat{A}}{\partial \ohat{f}_{i}} \, (\mass^{-1})_{ij} \, \ML_{jk}(\ohat{f})  \, (\mass^{-1})_{k\ell}\dfrac{\partial \hat{B}}{\partial \ohat{f}_{\ell}}.
\end{align}
The symmetric matrix $\ML(\ohat{f})$ is given by
\begin{multline}\label{eq:discrete-landau-matrix}
\ML_{ij}(\ohat{f})=-\dfrac{1}{2}\int \limits_{\Omega} \int \limits_{\Omega} \bigg( \dfrac{\partial \phy_{i} (v')}{\partial v'}-\dfrac{\partial \phy_{i} (v'')}{\partial v''} \bigg) \\
 \cdot  M(\f_h(v')) \, U(v'; v'') \, M(\f_h(v'')) 
\cdot\, \bigg( \dfrac{\partial \phy_{j} (v')}{\partial v'}-\dfrac{\partial \phy_{j} (v'')}{\partial v''} \bigg) \dv' \dv'' .
\end{multline}
It remains negative semi-definite, with a sufficient but not
necessary condition provided by the positive semi-definiteness of
$M(\f_h)$. This property will be used later on to prove the second
law of thermodynamics.
The action of the discrete metric bracket on two functionals,
$\hat{A} (\ohat{f}) = \fa [\f_{h}]$ and
$\hat{B} (\ohat{f}) = \fb [\f_{h}]$, can now be expressed as
\begin{align}\label{eq:dfd_discrete_bracket_gradient_form}
( \hat{A} , \hat{B} )_{h} = \nabla \hat{A}^{T} \MG (\ohat{f}) \, \nabla \hat{B} ,
\end{align}
where the gradient $\nabla$ is to be taken with respect to the degrees
of freedom $\ohat{f}$ and the matrix operator $\MG$ is given by
\begin{align}\label{eq:metric_bracket_discrete_matrix_kernel}
\MG (\ohat{f}) = \mass^{-1} \, \ML(\ohat{f}) \, \mass^{-1}.
\end{align}
Inserting $\hat{A}(\ohat{f})=\ohat{f}$ and $\hat{B}(\ohat{f})=\hat{F}$ into the bracket, the equations of motion for $\ohat{f}$ become
\begin{align}\label{eq:dfd_eqs_of_motion_gradient_form}
\dfrac{d\ohat{f}}{dt} = \MG (\ohat{f}) \, \nabla \hat{F} = ( \ohat{f} , \hat{F} )_{h},
\end{align}
in direct analogy with the continuous case.

\subsection{Semi-discrete Conservation Laws}
The mass, $\fm$, momentum, $\mom$, and kinetic energy, $\nrg$, carried by the
distribution function $\f$ are defined in terms of the functionals
\begin{align}
\fm[\f]  = \int \limits_{\Omega} \, \f(v) \dv, \qquad
\mom[\f] = \int \limits_{\Omega} \, v \, \f(v) \dv, \qquad
\nrg[\f] = \dfrac{m}{2} \int \limits_{\Omega} \, \vert v\vert^2 \, \f(v) \dv.
\end{align}
Analogously, the mass, momentum and kinetic energy carried by the discrete
distribution function $\f_{h}$ are
\begin{subequations}\label{eq:semi-discrete-casimirs}
\begin{align}
\hat{M} (\ohat{f}) &\equiv \fm [f_h] = \sum \limits_{i=1}^{N} \ohat{f}_i \int \limits_{\Omega} \, \phy_i(v) \dv, \\
\hat{P} (\ohat{f}) &\equiv \mom[f_h] = \sum \limits_{i=1}^{N} \ohat{f}_i \int \limits_{\Omega} \, v \, \phy_i(v) \dv, \\
\hat{E} (\ohat{f}) &\equiv \nrg[f_h] = \dfrac{m}{2} \sum \limits_{i=1}^{N} \ohat{f}_i \int \limits_{\Omega} \, \vert v\vert^2 \, \phy_i(v) \dv.
\end{align}
\end{subequations}
If $Q_h(\Omega)$ is such that $\{ 1, v, v^2 \} \in Q_h(\Omega)$, which is the case e.g. for quadratic Lagrange finite elements, we can find coefficients $\ohat{v}$ and $\ohat{\varepsilon}$, so that
\begin{align}
\sum \limits_{i=1}^{N} \ohat{v}_i \phy_i(v) = v , \qquad
\sum \limits_{i=1}^{N} \ohat{\varepsilon}_i \phy_i(v) = v^2 ,
\end{align}
in terms of which the discrete quantities~\eqref{eq:semi-discrete-casimirs} become
\begin{align}\label{eq:semi-discrete-casimirs-2}
\hat{M} (\ohat{f})
= \unity^T \mass \, \ohat{f},\qquad
\hat{P} (\ohat{f})
= \ohat{v}^T \mass \, \ohat{f},\qquad
\hat{E} (\ohat{f})
= \dfrac{m}{2} \, \ohat{\veps}^T \mass \, \ohat{f} ,
\end{align}
where $\unity$ denotes the vector in $\rsp^{N}$ with all elements being $1$.
The mass, momentum and kinetic energy carried by $\ohat{f}$ are then
Casimirs, and thus invariants of motion, of the discrete metric
bracket~\eqref{eq:dfd_discrete_bracket_gradient_form} since, for
arbitrary $\hat{A}$, we have
\begin{subequations}
\begin{align}
( \hat{M} , \hat{A} )_{h} = \unity^T \mass G (\ohat{f}) \, \nabla \hat{F}  = \unity^{T} \ML(\ohat{f}) \mass^{-1} \nabla \hat{A} = 0,\\
( \hat{P} , \hat{A} )_{h} = \ohat{v}^T \mass G (\ohat{f}) \, \nabla \hat{F}  = \ohat{v}^{T} \ML(\ohat{f}) \mass^{-1} \nabla \hat{A} = 0,\\
( \hat{E} , \hat{A} )_{h} = \ohat{\veps}^T \mass G (\ohat{f}) \, \nabla \hat{F} = \ohat{\veps}^{T} \ML(\ohat{f}) \mass^{-1} \nabla \hat{A} = 0.
\end{align}
\end{subequations}
This follows directly from the properties of the Landau matrix $\ML$, namely
\begin{subequations}
\begin{align}
\unity^{T} \ML(\ohat{f}) &\equiv 0 , \textrm{ for any } \ohat{f}, \\
\ohat{v}^{T} \ML(\ohat{f}) &\equiv 0 , \textrm{ for any } \ohat{f}, \\
\ohat{\veps}^{T} \ML(\ohat{f})& \equiv 0 , \textrm{ for any } \ohat{f}.
\end{align}
\end{subequations}
Note that here we only consider the kinetic energy $\mcal{E}$ instead of the whole Hamiltonian $\ham$. However, the metric bracket does not act on the electric and magnetic field, therefore the discussion of the current and the next section does not change when replacing one with the other.

\subsection{Semi-discrete H-Theorem}
\label{sec:semi-discrete-h-theorem}

We now want to reproduce the discussion of Section~\ref{sec:general} on the discrete level, focusing however on the metric case.
The time evolution of the discrete entropy $\hat{S} (\ohat{f}) \equiv \fen[f_h]$ is given by
\begin{align}
\dfrac{d \hat{S}}{dt} = ( \hat{S}, \hat{F} )_{h} = - ( \hat{S} , \hat{S} )_{h} \ge 0,
\end{align}
where $\hat{F} = \hat{E} (\ohat{f}) - \hat{S} (\ohat{f})$ with $\hat{E} (\ohat{f}) \equiv \nrg[f_h]$, and we used the fact that the kinetic energy $\hat{E}$ is in the nullspace of the discrete bracket.
The last inequality follows from the fact that $\ML(\ohat{f})$, and thus the discrete bracket, is negative semi-definite as long as $M(\f_h)$ remains positive semi-definite. Thus the discrete entropy evolves monotonically in time under the action of the discrete metric bracket.

As in the continuous case, the equilibrium state is determined by the energy-Casimir method, that is by requiring
\begin{align}
\de \hat{S} (\ohat{f}_{eq}) + \de \sum \limits_{i} \lambda_{i} \hat{C}_{i} (\ohat{f}_{eq}) = 0 ,
\end{align}
or, specifically, for the discrete Casimir invariants as given in~\eqref{eq:semi-discrete-casimirs-2},
\begin{align}
\bigg(
\de \hat{S}
+ \lambda_{\hat{M}} \, \de (\unity^T \mass \, \ohat{f})
+ \lambda_{\hat{P}} \, \de (\ohat{v}^T \mass \, \ohat{f})
+ \lambda_{\hat{E}} \, \de \bigg( \dfrac{m}{2} \, \ohat{\veps}^T \mass \, \ohat{f} \bigg)
\bigg) \bigg\vert_{\ohat{f} = \ohat{f}_{eq}} = 0 .
\end{align}
With the discrete entropy functional given by
\begin{align}
\hat{S} (\ohat{f}) \equiv \fen [f_{h}] = \int \limits_{\Omega} s \bigg( \sum \limits_{j} \ohat{f}_{j} \phi_{j} (v) \bigg) \dv ,
\end{align}
we obtain the following condition for the discrete equilibrium state $\ohat{f}_{eq}$,
\begin{align}
\sum \limits_{i} \bigg[ \nabla_{i} \hat{S} (\ohat{f}_{eq})
+ \sum \limits_{j} \bigg(
     \lambda_{\hat{M}} \, \mass_{ji}
   + \lambda_{\hat{P}} \, \ohat{v}_{j} \mass_{ji}
   + \dfrac{m}{2} \lambda_{\hat{E}} \, \ohat{\veps}_{j} \mass_{ji}
\bigg) \bigg] \, \delta \ohat{f}_{eq,i} = 0 ,
\end{align}
where the gradient of the discrete entropy function,
\begin{align}
\nabla_{i} \hat{S} (\ohat{f}) = \int \limits_{\Omega} \phi_{i} (v) \, s_{\f} \bigg( \sum \limits_{j} \ohat{f}_{eq,j} \phi_{j} (v) \bigg) \dv ,
\end{align}
corresponds to a projection of $s_{\f} (\f_{h})$ onto the finite-dimensional space $Q_h(\Omega)$.
In order to obtain the discrete equilibrium state, we have to solve
\begin{align}
\nabla \hat{S} (\ohat{f}_{eq})
+ \lambda_{\hat{M}}  \, \mass \unity
+ \lambda_{\hat{P}} \, \mass \ohat{v}
+ \dfrac{m}{2} \, \lambda_{\hat{E}} \, \mass \ohat{\veps} = 0 ,
\end{align}
and by matching the mass, momentum and kinetic energy of the equilibrium and the initial condition, we can determine the multipliers $\lambda_{\hat{M}}$, $\lambda_{\hat{P}}$ and $\lambda_{\hat{E}}$, respectively.
Uniqueness of the discrete equilibrium state $\ohat{f}_{eq}$ follows from the same convexity argument as in the continuous case.

\section{Temporal Discretization}

For the temporal discretization, we consider so-called discrete
gradients. We briefly summarize the essence of these methods and then
prove the conservation laws for the fully discrete system. Finally, we
rearrange the resulting nonlinear equations into a sparse form that is
suitable for standard algorithms addressing large nonlinear systems.

\subsection{Discrete Gradients}

Discrete gradients~\cite{Quispel:1996, McLachlan:1999, QuispelMcLaren:2008, MansfieldQuispel:2009} constitute discrete analogues of the gradient of a function,
and generalizations thereof. These methods can be applied to any
system of ordinary differential equations in the form
\begin{align}
\dfrac{du}{dt} = \MS(u) \, \nabla F(u) ,
\end{align}
where $\MS(u)$ can be an anti-symmetric matrix for conservative systems,
a symmetric matrix for dissipative systems, or a combination thereof
for metriplectic systems.  This obviously resembles the structure
of~\eqref{eq:dfd_eqs_of_motion_gradient_form}.  This system of ODEs is
discretized by
\begin{align}\label{eq:discrete_gradient}
\dfrac{u_{n+1} - u_{n}}{\Delta t} = \bar{\MS} (u_{n}, u_{n+1}) \, \bar{\nabla} F (u_{n}, u_{n+1}) ,
\end{align}
where $\bar{\MS} (u_{n}, u_{n+1})$ is any symmetric or anti-symmetric
matrix that approaches $\MS(u)$ in the limit of
$u_{n+1} \rightarrow u_{n}$ and $\Delta t \rightarrow 0$,
so that $\bar{\MS} (u_{n}, u_{n}) = \MS(u_{n})$, and
$\bar{\nabla} F (u_{n}, u_{n+1})$ is a discrete gradient.  Given a
differentiable function $F : \rsp^{m} \rightarrow \rsp$, a discrete
gradient $\bar{\nabla} F (u_{n}, u_{n+1})$ is a vector valued
continuous function of $(u_{n}, u_{n+1})$, satisfying
\begin{align}
\label{eq:discrete-gradient-property}
\begin{split}
(u_{n+1} - u_{n}) \cdot \bar{\nabla} F (u_{n}, u_{n+1}) &= F (u_{n+1}) - F (u_{n}) , \\
\bar{\nabla} F (u_{n}, u_{n}) &= \nabla F (u_{n}) .
\end{split}
\end{align}
Several such discrete gradients are known. One may consider for
example the midpoint discrete gradient by \citet{Gonzalez:1996},
\begin{multline}
\bar{\nabla} F (u_{n}, u_{n+1})
= \nabla F (u_{n+1/2}) \\
+ ( u_{n+1} - u_{n} ) \, \dfrac{F (u_{n+1}) - F (u_{n}) - ( u_{n+1} - u_{n} ) \cdot \nabla F (u_{n+1/2}) }{\abs{ u_{n+1} - u_{n} }^{2}} , 
\end{multline}
with $u_{n+1/2} = \tfrac{1}{2} (u_{n} + u_{n+1})$, or the average
discrete gradient by \citet{HartenLaxLeer:1983},
\begin{align}
\bar{\nabla} F (u_{n}, u_{n+1}) &= \int \limits_{0}^{1} \nabla F \big( (1-\xi) u_{n} + \xi u_{n+1} \big) \, d\xi .
\end{align}
Higher-order time integrators can be constructed following~\citet{CohenHairer:2011}.

\subsection{Discrete Conservation Laws}
Using the concept of discrete gradients for the temporal discretization allows us to prove momentum and energy conservation as well as the correct monotonic behavior of the entropy.
For definiteness, we consider a midpoint discretization for $\bar{\MS}$, i.e.,
\begin{align}
\bar{\MS} (\ohat{f}_{n}, \ohat{f}_{n+1}) &= \MS (\ohat{f}_{n+1/2}) , &
& \text{where} &
\ohat{f}_{n+1/2} &= \dfrac{\ohat{f}_{n} + \ohat{f}_{n+1}}{2} ,
\end{align}
but let us note that other choices work just as well.
For the discrete metric system, the matrix $\MS$ is given by the matrix operator $\MG$ in~\eqref{eq:metric_bracket_discrete_matrix_kernel}, that is
\begin{align}
\MS (\ohat{f}) =\mass^{-1} \, \ML (\ohat{f}) \, \mass^{-1} ,
\end{align}
and $F$ corresponds to the discrete free energy $\hat{F} (\ohat{f}) = \hat{E} (\ohat{f}) - \hat{S} (\ohat{f})$.
Assuming a discrete gradient time stepping algorithm~\eqref{eq:discrete_gradient}, the difference of the energy at consecutive points in time is obtained from
\begin{align}
\ohat{\varepsilon}^T \mass (\ohat{f}_{n+1} - \ohat{f}_{n})
\nonumber
&= \Delta t \, \ohat{\varepsilon}^T \, \mass \, \bar{\MS} (\ohat{f}_{n}, \ohat{f}_{n+1}) \, \hat{F} (\ohat{f}_{n}, \ohat{f}_{n+1}) \\
&= \Delta t \, \ohat{\varepsilon}^T \, \mass \, \mass^{-1} \, \ML (\ohat{f}_{n+1/2}) \, \mass^{-1} \, \bar{\nabla} \hat{F} (\ohat{f}_{n}, \ohat{f}_{n+1}) .
\end{align}
The right-hand side vanishes exactly since $\ohat{\varepsilon}^T \mass \mass^{-1} \ML (\ohat{f}_{n+1/2}) = \ohat{\varepsilon}^T \ML (\ohat{f}_{n+1/2})=0$. 
Therefore we have that $\ohat{\varepsilon}^T \mass \ohat{f}_{n+1} = \ohat{\varepsilon}^T \mass \ohat{f}_{n}$, stating that the energy at time $t_{n+1}$ equals the energy at time $t_{n}$.
The momentum and density conservation are proved in full analogy.

\subsection{Discrete Entropy Production}

The monotonic increase of entropy can be shown as follows.
The difference of the entropy at two consecutive points in time is
\begin{align}
\hat{S}(\ohat{f}_{n+1}) - \hat{S}(\ohat{f}_{n}) &= \big( \hat{E}(\ohat{f}_{n+1}) - \hat{F}(\ohat{f}_{n+1}) \big) - \big( \hat{E}(\ohat{f}_{n}) - \hat{F}(\ohat{f}_{n}) \big) ,
\end{align}
where the discrete energy function $\hat{E}$ is conserved, so that $\hat{E}(\ohat{f}_{n+1}) = \hat{E}(\ohat{f}_{n})$.
Using property~\eqref{eq:discrete-gradient-property} of the discrete gradient, we have
\begin{align}
\hat{S}(\ohat{f}_{n+1})-\hat{S}(\ohat{f}_{n})
\nonumber
&= - \big( \hat{F}(\ohat{f}_{n+1})-\hat{F}(\ohat{f}_{n}) \big) \\
\nonumber
&= - \Delta t\, \bar{\nabla} \hat{F}^T (\ohat{f}_{n}, \ohat{f}_{n+1}) \, \mass^{-1} \, \ML (\ohat{f}_{n+1/2}) \, \mass^{-1}\, \bar{\nabla} \hat{F} (\ohat{f}_{n}, \ohat{f}_{n+1}) \\
&\ge 0 ,
\end{align}
as both $\ML$ is symmetric negative semi-definite matrix and $\mass$ is a symmetric matrix.

\subsection{Sparsification for Nonlinear Solver}
The introduction of discrete gradient methods for temporal discretization results in a nonlinear system of equations for the degrees of freedom $\ohat{f}_{n+1}$. Such systems are typically solved with Jacobian-free quasi-Newton methods accompanied with Krylov-subspace iteration to solve the resulting linear system at each quasi-Newton step. 

The time and space discrete system involves the Landau matrix, $\ML(\ohat{f})$, and the inverse of the mass matrix, $\mass^{-1}$, which are both dense matrices for typical finite element discretizations. Storing such matrices could quickly require extensive memory resources, rendering the matrix-vector products, that are required for the Krylov-subspace construction, expensive and slow. To facilitate an efficient implementation of the algorithm, a formulation involving only sparse matrices is thus desirable.

This can be accomplished by first multiplying the discrete system for $\ohat{f}_{n+1}$ with the mass matrix from left according to
\begin{align}
\mass\, (\ohat{f}_{n+1} - \ohat{f}_{n})
= \Delta t \, \ML (\ohat{f}_{n+1/2}) \, \mass^{-1} \, \bar{\nabla} \hat{F} (\ohat{f}_{n}, \ohat{f}_{n+1}),
\end{align}
and then rearranging the right-hand side of the above equation into 
\begin{align}
&\ML (\ohat{f}_{n+1/2}) \, \mass^{-1} \, \bar{\nabla} \hat{F} (\ohat{f}_{n}, \ohat{f}_{n+1})=\MC(\ohat{f}_n,\ohat{f}_{n+1})\ohat{f}_{n+1/2} ,
\end{align}
so that the final system of equations for $\ohat{f}_{n+1}$ reads
\begin{align}\label{eq:sparse-system}
\mass\, (\ohat{f}_{n+1} - \ohat{f}_{n})=\Delta t\, \MC(\ohat{f}_n,\ohat{f}_{n+1})\, \ohat{f}_{n+1/2}.
\end{align}
The matrix $\MC$ is sparse and given by
\begin{align}
\MC(\ohat{f}_n,\ohat{f}_{n+1})=&\int \limits_{\Omega}  \dfrac{\partial \phy_{i} (v')}{\partial v'} \cdot D(\ohat{f}_n,\ohat{f}_{n+1})  
\cdot \left(\dfrac{\partial \phy_{k} (v')}{\partial v'}(\mass^{-1})_{k\ell} \, \bar{\nabla}_{\ell} \hat{F} (\ohat{f}_{n}, \ohat{f}_{n+1})\right)\phy_{j} (v')\dv'
\nonumber\\
&-\int \limits_{\Omega}\dfrac{\partial \phy_{i} (v')}{\partial v'} \cdot K(\ohat{f}_n,\ohat{f}_{n+1})\,  \phy_{j} (v') \,\dv',
\end{align}
where the tensor $D$ and the vector $K$, required to construct $\MC$, are functions of the velocity coordinate and determined according to
\begin{align}
&D(\ohat{f}_n,\ohat{f}_{n+1})=\int \limits_{\Omega} U (v', v'') \,\f_{h,n+1/2}(v'')\dv'',\\
&K(\ohat{f}_n,\ohat{f}_{n+1})=\int \limits_{\Omega}  U (v', v'') \,\f_{h,n+1/2}(v'')\cdot \left(\dfrac{\partial \phy_{k} (v'')}{\partial v''}(\mass^{-1})_{k\ell} \, \bar{\nabla}_{\ell} \hat{F} (\ohat{f}_{n}, \ohat{f}_{n+1})\right) \dv''.
\end{align}

In its new form, \eqref{eq:sparse-system}, all matrices in the nonlinear system of equations for the time advance are sparse. A numerically estimated Jacobian for any quasi-Newton method will be sparse, and the matrix-vector products for constructing the Krylov subspace become feasible. The sparse form has the same structure as in \cite{Hirvijoki:2017ei} were the solvability was demonstrated numerically, although the formulation in general was not based on the metriplectic formulation.

\section{Summary and Outlook}
In this work, a general framework for metriplectic Galerkin discretizations
of the Landau collision integral was presented. The discretization proceeded in two steps.
First, a semi-discretization of the metriplectic bracket for the collision integral was obtained, which preserves
the Casimirs of the metric bracket and guarantees entropy to be a monotonic function of time. The semi-discrete system is thus a finite-dimensional metriplectic system. Then, the system was discretized in time by
employing the discrete gradients methods, still retaining exact conservation of Casimirs and the monotonic behavior of entropy. Therefore the resulting method corresponds to one of the rare instances of a genuine metriplectic integrator. 

One of the advantages of our approach is that the conservation laws end entropy behavior are not
manufactured or forced into the scheme ``by hand'' but follow automatically from preserving the underlying metriplectic structure of the Landau collision integral. Another advantage is the generality of the spatial discretization, allowing one to employ the full power of finite element discretizations and non-structured meshes.

While numerical demonstrations are left for future work, we discuss in detail how to convert the initially dense nonlinear system of equations into a sparse form that is suitable for Jacobian-free Newton-Krylov methods. The numerical implementation and the metriplectic discretization of the full Vlasov--Maxwell--Landau system are currently under investigation and will be addressed in future publications.

Due to the flexibility of both the metriplectic formalism and our discretization approach, many potential generalizations are possible. Obviously, other collision operators that fit into the metriplectic framework are immediately treatable.
The fact that for any desired equilibrium state there is some freedom in the construction of the bracket might also turn out to be useful, e.g., for the design of more refined numerical algorithms, that use this freedom to improve efficiency or efficacy.
Further, it might be possible to construct discrete metric brackets that relax the system to more complicated equilibrium states, which might be useful for constraint relaxation problems or for the dissipation of small scale structures by choosing an evolving equilibrium that corresponds, e.g., to a coarse grained solution of the current distribution function.

\section*{Acknowledgments}

\noindent Useful discussions with Omar Maj and Philip Morrison are gratefully acknowledged.
MK has received funding from the European Union's Horizon 2020 research and innovation programme under the Marie Sklodowska--Curie grant agreement No 708124. EH was supported by the U.S. Department of Energy Contract No. DE-AC02-09-CH11466. The views and opinions expressed herein do not necessarily reflect those of the European Commission or the U.S. Department of Energy.

\bibliographystyle{unsrtnat}
\bibliography{mi_landau}

\end{document}